\documentclass[12pt]{article}

\usepackage{amsfonts}

\title{On Sets of Integers where Each Pair Sums to a Square }

\author{
Allan J. MacLeod\\Dept. of Mathematics and Statistics,\\
University of the West of Scotland,\\High St.,  Paisley,\\Scotland.  PA1 2BE\\
(e-mail: allan.macleod@uws.ac.uk) }

\date{}

\begin{document}

\maketitle

\begin{abstract}
{\noindent We discuss the problem of finding distinct integer sets
$\{x_1,x_2,\ldots,x_n\}$ where each sum $x_i+x_j, \, i \ne j$ is
a square, and $n \le 7$. We confirm minimal results of Lagrange
and Nicolas for $n=5$ and for the related problem with triples. We
provide new solution sets for $n=6$ to add to the single known
set. This provides new information for problem D15 in Guy's {\it
Unsolved Problems in Number Theory} }

\vspace{0.5cm}

\end{abstract}

\section{Introduction}
Let $\{x_1,x_2,\ldots,x_n\}$ be a set of $n$ distinct non-zero
integers. Erdos \cite{erd} and, independently, Moser, in
\cite{sier}, asked for examples of sets where $x_i+x_j$ is a
square for all possible pairs of subscripts with $i \ne j$. This
is considered in section D15 of Guy's well-known book \cite {guy}.

For $n=2$, there is only $1$ pair and so the problem is easy.
Pick any square $p^2$, then $\{x,p^2-x\}$ is a suitable set.

For $n=3$, we have a system of $3$ equations in $3$ unknowns,
\begin{equation}
\left( \begin{array}{ccc}1&1&0\\1&0&1\\0&1&1 \end{array} \right )
\left( \begin{array}{c}x_1\\x_2\\x_3 \end{array} \right ) = \left
( \begin{array}{c}p^2\\q^2\\r^2 \end{array} \right )
\end{equation}

This is easily solved to give
\begin{equation}
\left( \begin{array}{c}x_1\\x_2\\x_3 \end{array} \right ) =
\frac{1}{2} \left(
\begin{array}{c}p^2+q^2-r^2\\p^2-q^2+r^2\\-p^2+q^2+r^2
\end{array} \right )
\end{equation}

Since squares are either congruent to $0$ or $1$ modulo $4$, we
cannot have more than one $x_i$ of the form $4k+1$ and more than
one of the form $4k+3$. Thus, at least one value must be even,
either of the form $4m$ or $4m+2$. If we had both types of odd
numbers, then the sum of the even value plus the odds would give
at least one value of the form $4i+3$ which cannot be a square.
Thus, the set can consist of at most one odd value, with the
remainder even. Clearly, at most one value can be negative. This
is true for all larger sets.

If we choose $3$ even squares, we can get a set of even $x_i$,
for example, $(2^2,4^2,8^2)$ gives the set $\{-22,26,38\}$,
whilst, if we choose $2$ odd squares and $1$ even square we get a
set with one odd element, for example, $(1^2,2^2,3^2)$ gives
$\{-2,3,6\}$. Other choices lead to numbers with denominator $2$,
and if we multiply all $x_i$ by $4$ we preserve the paired
squareness.

For $n=4$, we follow the description given by Lagrange \cite{lag1}
and Nicolas \cite{nic}. We have $6$ equations,
\begin{equation}
\begin{array}{ccc}
x_1+x_2=p^2& \hspace{2cm} & x_3+x_4=t^2 \\
x_1+x_3=q^2& \hspace{2cm} & x_2+x_4=u^2 \\
x_2+x_3=r^2& \hspace{2cm} & x_1+x_4=v^2
\end{array}
\end{equation}
and so $S=x_1+x_2+x_3+x_4=p^2+t^2=q^2+u^2=r^2+v^2$.

Thus we, first, look for numbers $S$ which have at least $3$
different representations as the sum of $2$ squares. Let the
prime decomposition of $S$ be
\begin{displaymath}
S= 2^i \; p_1^{n_1}p_2^{n_2}\ldots p_k^{n_k} \;
q_1^{m_1}q_2^{m_2} \ldots q_l^{m_l}
\end{displaymath}
where the $p_i$ primes are $\equiv 3 \bmod 4$ and the $q_i$ primes
are $\equiv 1 \bmod 4$.

Then, if one of the $n_i$ values is odd, the number of
representations is $0$. If all $n_i$ are even, then the number of
representations is
\begin{displaymath}
\frac{(m_1+1)(m_2+1) \ldots (m_l+1) + \delta}{2}
\end{displaymath}
where $\delta =1$ if all $m_i$ are odd and $\delta=0$ otherwise.

Having found a suitable $S$, we consider all possible groups of
$3$ different representations. The left hand set of equations in
(3) is solved as in equation (2), and $x_4$ computed from the
right-hand-sides. By changing the order of $r$ and $v$ we derive
a second solution. All other permutations give one of these two
basic solutions.

\begin{center}
TABLE 1\\Smallest sets of $4$ elements\\
\begin{tabular}{rrrr}
$\;$&$\;$&$\;$&$\;$ \\
$x_1$&$x_2$&$x_3 \,$&$x_4 \,$\\
-40&65&104&296 \\
-94&95&130&194 \\
-88&88&137&488 \\
-94&98&263&578 \\
-190&239&290&386
\end{tabular}
\end{center}

\begin{center}
TABLE 2\\Smallest sets of $4$ positive elements\\
\begin{tabular}{rrrr}
$\;$&$\;$&$\;$&$\;$ \\
$x_1$&$x_2 \;$&$x_3 \;$&$x_4 \; \,$\\
2&359&482&3362\\
8&1016&1288&3473\\
162&567&1282&4194\\
2&167&674&6722\\
98&863&1346&5378
\end{tabular}
\end{center}

All this is very easily programmed using the freely available
system Pari-GP, and can be run very quickly on any modern machine.
These calculations lead to the results in Tables $1$ and $2$,
which give the $5$ smallest general sets and strictly positive
sets, respectively.

To measure the size of solutions, we use the $l_1$ norm $\sum
|x_i|$ rather than $S=\sum x_i$, as used by Lagrange and Nicolas.
The different measure of size means that the smallest set given
by Lagrange and Nicolas is only the second smallest in Table 1.

\section{n=5 sets}
Suppose $\{x_1, x_2, x_3, x_4\}$ is a set such that all $6$ pairs
$x_i+x_j$, with $i \ne j$, are squares. We assume $x_1 < x_2 <
x_3 < x_4$, and look for a fifth element $x_5$ to make a set of
five elements.

We know $p$ and $q$ such that $x_1+x_4=p^2$ and $x_2+x_4=q^2$, and
look for $x_5$ with $x_1+x_5=w^2$ and $x_2+x_5=y^2$. Thus
$x_2-x_1=q^2-p^2=y^2-w^2=(y+w)(y-w)$. We loop over the divisors of
$q^2-p^2$, find suitable $y$ and $w$, and then find $x_5=w^2-x_1$,
checking if it satisfies the other conditions, namely $x_3+x_5$
and $x_4+x_5$ being square.

\begin{center}
TABLE 3\\Smallest sets of $5$ elements\\
\begin{tabular}{rrrrr}
$\;$&$\;$&$\;$&$\;$ \\
$x_1 \,$&$x_2 \;$&$x_3 \, \;$&$x_4 \; \;$&$x_5 \; \; \,$\\
 -4878& 4978& 6903& 12978& 31122\\
  -2158& 2258& 4967& 19058& 33842\\
  -5998& 7847& 9842& 11474& 30962\\
  -878& 882& 7767& 12114& 48402\\
  -1417& 1586& 5138& 18578& 45938
\end{tabular}
\end{center}

Tables $3$ and $4$ give the smallest sets (in the $l_1$ norm
sense) of general integers and positive integers, respectively.
Table 3 agrees with the results in Lagrange apart from the second
and third sets being swapped, and the second smallest in Table 4
confirms Lagrange's speculation, since his table of results only
contains one all-positive set, the first one in Table 4.

\begin{center}
TABLE 4\\Smallest sets of $5$ positive elements\\
\begin{tabular}{rrrrr}
$\;$&$\;$&$\;$&$\;$ \\
$x_1 \, \,$&$x_2 \; \,$&$x_3 \; \; \,$&$x_4 \; \; \;$&$x_5 \; \; \; \,$\\
 7442& 28658& 148583& 177458& 763442\\
  32018& 104882& 188882& 559343& 956018\\
  9122& 104447& 208034& 348482& 1295042\\
  23458& 82818& 127863& 228546& 2149218\\
  30818& 56207& 322018& 910082& 1946018
\end{tabular}
\end{center}

Now since the binomial coefficients $\left( \begin{array}{c}5\\2
\end{array} \right) = \left( \begin{array}{c}5\\3
\end{array} \right)$, there are the same number of subsets with
$3$ elements as with $2$ elements, and we can relate them easily.

If we let $S= \sum x_i$, define the values $z_i=S/3-x_i,
i=1,\ldots,5$, then it is easy to see that the sum of any $3$
non-repeated $z_i$ values equals the sum of two of the $x_i$
values, using the remaining subscripts eg. $z_1+z_3+z_4=x_2+x_5$.
Thus $\{z_1, z_2, z_3, z_4, z_5 \}$ forms a set of $5$ elements
where every triple sums to a square. If we have values of $z_i$
which are rational we can scale by $9$ to provide integer sets.

\begin{center}
TABLE 5\\Smallest sets with square triplets and $5$ positive elements\\
\begin{tabular}{rrrrr}
$\;$&$\;$&$\;$&$\;$ \\
$x_1 \; \; \,$&$x_2 \; \; \,$&$x_3 \; \; \,$&$x_4 \; \; \;$&$x_5 \; \; \; \,$\\
92763& 4914963& 7559299& 9945963& 16308963\\
1039923& 2292723& 5649363& 10128915& 21847678\\
695883& 2655435& 18466923& 40327563& 62161518\\
33843& 22986003& 75168435& 123438558& 167502963\\
22906587& 36372270& 71091867& 114486267& 211586907
\end{tabular}
\end{center}

This problem was considered by Gill \cite{gill} in his wonderful
book, and then in a more modern context by Wagon \cite{wag}, who
published a positive solution where each $z_i$ has at least $20$
digits. The computer search which produced Tables $3$ and $4$ can
be easily adapted to search for positive sets - small sets with
negatives arise easily from small sets in Table 3. These triplet
results are in Table 5.

These values are significantly smaller than that found by Wagon. In
fact, the first set is given by Lagrange, right at the end of
\cite{lag2}, though it is unclear which of Lagrange and Nicolas
discovered the solution. The present results agree with the claim
made that the first set is the smallest positive set giving square
triples. This answers the question in Section D15 of Guy \cite{guy}.
It is not obvious why this solution has lain unnoticed for so long.
My $4$ years of school French are enough to give a reasonable
understanding of the papers.

Comparing Tables $4$ and $5$ we see that the values for positive
sets with square triples are significantly higher than for sets
with square pairs. A possible reason comes from the following
analysis.

If we assume that $x_1 < x_2 < x_3 < x_4 < x_5$, then the
smallest of the $z_i$ values is $z_5$. So what we want to know is
$\Pr(z_5>0$). But this is equal to
\begin{description}
\item  $\Pr(\, S > 3x_5 \,)=\Pr(\, x_1+x_2+x_3+x_4 > 2x_5 \,)$
\item $=\Pr(\, \,3x_1 + x_2 + x_3 + x_4 > 2( x_1 + x_5) \, \, )$
\item $=\Pr(\, \,(x_1+x_2) \,+ \,(x_1+x_3) \, + \,(x_1+x_4) > 2(x_1+x_5)\, \,)$
\item $=\Pr(\, p^2 \, + \, q^2 \, + \, r^2 \, > \, 2s^2 \,)$
\end{description}
for some positive integers $p,q,r,s$.

Now $p < q < r < s$, so what we want is $\Pr(\, y_1^2 \, + \,
y_2^2 \, + \, y_3^2 \, > \, 2 \,)$ where $0 < y_1 < y_2 < y_3 <
1$. Unfortunately, the distribution of the $y_i$ is unclear to
me, so to proceed I assume they are uniformly distrubuted in
$[0,1]$.

Given a random point in the unit cube, there are $6$ possible
orderings of the coordinates, so the probability with $y_1 < y_2
< y_3$ is $1/6$ of the probability for a general point.

The sphere $\, y_1^2 + y_2^2 + y_3^2 = 2\,$ cuts the unit cube at
the obvious points $(1,1,0), (1,0,1), (0,1,1)$ and lies both
inside and outside the cube, for example $(1,1,1)$ is outside.
Thus the required probability is just the volume within the cube
but outside the sphere. The volume inside both cube and sphere is
\begin{displaymath}
\int_0^1 \, \, \int _0^1 \, \, \int_0^{\min\{1,\sqrt{2-x^2-y^2}\}}
\, \, dz \, dy \, dx
\end{displaymath}
which can be written as
\begin{displaymath}
\int_0^1 \int_0^{\sqrt{1-x^2}} \, dy \, dx \, \, + \, \, \int_0^1
\, \int_{\sqrt{1-x^2}}^1 \, \, \sqrt{2-x^2-y^2} \, dy \, dx
\end{displaymath}

The first integral is easy and gives the answer $\pi / 4$. The
second is more difficult but possible, and was checked with a
symbolic package giving $\pi(1-2\sqrt{2}/3)$. Thus the required
probability (taking into account the $1/6$ term) is
$(\pi(8\sqrt{2}-15)+12)/72 \approx 1/172$. This gives a
reasonable explanation of the rarity of positive triples.

\section{n=6 Sets}
We can apply the ideas of the last section to finding a value
$x_6$ to give a set of $6$ elements. There are $15$ possible
pairings. Lagrange calls a set satisfying only $14$ a
pseudo-solution. These are fairly easy to find and a large table
is given in \cite{lag2}.

A special observation about such pseudo-solutions is that many
have $S=\sum x_i$ being an integer square. Since $S$ consists of
three pairs which should each sum to a square, this leads to a
consideration of representations of the form $ p^2 + q^2 + r^2 =
s^2$. Lagrange uses the identity
\begin{displaymath}
(t^2+u^2+v^2+w^2)^2=(t^2+u^2-v^2-w^2)^2+4(t\;w-uv)^2+4(tv+uw)^2
\end{displaymath}
and a very inventive magic-square approach to find a
representation which easily gives $13$ out the the $15$ pairs
adding to a square. The remaining $2$ pairs require two binary
quartics to be made square. Each quartic is of the form
\begin{displaymath}
a \; G^4 + b \; G^3 H + c \; G^2 H^2 - b \; G H^3 + a \; H^4 =
F^2,
\end{displaymath}
though, obviously, with different $(a,b,c)$ values.

Both quartics have many solutions and are thus birationally
equivalent to different elliptic curves, which, on investigation,
have large ranks of the order $3, 4, \ldots$. Thus, there will be
lots of solutions with the possibility of a common solution,
leading to a set satisfying all $15$ identities. Lagrange
completed the square, in the style of Euler, of one of the
quartics, and used the resulting necessary condition to simplify
the other quartic which he was again able to complete the square.
This led to the solution set
\begin{displaymath}
\{-15863902, 17798783, 21126338, 49064546, 82221218, 447422978\}
\end{displaymath}

In the current project, we searched the first quartic for values
of $(G,H)$ satisfying the first quartic and simply tested whether
they satisfied the second. Some solutions give sets where values
are repeated, which we discarded. The style of quartic allows us
to only consider positive $G,H$ values, and other coding tricks
give a reasonably efficient code.

The program finds the solution of Lagrange very quickly and after
some more computing, we find the following further solutions
\begin{displaymath}
\{-1126475417645550, 1274834760326775, 2236102802190450,
\end{displaymath}
\begin{displaymath}
3990082124377234, 5054722548678450, 82588671432234450\},
\end{displaymath}

\begin{displaymath}
\{-5098887366661368, 6849471293061768, 8279927229562632,
\end{displaymath}
\begin{displaymath}
21501934179045768, 78740349884517393, 340192944008301832\}
\end{displaymath}
and
\begin{displaymath}
\{-1177836637755448, 2476350655765448, 2723928921099848,
\end{displaymath}
\begin{displaymath}
5744011161331073, 7858945782510152, 33438182171924552\}
\end{displaymath}

It is perfectly possible that we have missed some smaller
solutions. These results support Lagrange's conjecture that there
are an infinite number of different solution sets for $n=6$.

\section{$n=7$ Sets}
Given the small number of solutions found for the $n=6$ problem
it would have been very surprising if we had found a solution for
the $n=7$ problem, which would require $21$ pairs to be square.

Thus, we have concentrated in looking for sets which maximise the
number of pairs. We first looked at the $4$ solutions for the
$n=6$ problems and tried to extend them, using the method
described before.

We found that, if we add $15945698$ to the set found by Lagrange,
we have $18$ square pairs. None of the other $3$ solutions give
more than $17$ square pairs in an extension.

Next, we used the Lagrange parameterisation which gives $14$
square pairs and tried to add a seventh integer. This led to two
$18$ squares solutions
\begin{displaymath}
\{-256711392, 599109408, 741988233, 3222602992,
\end{displaymath}
\begin{displaymath}
 5845768992,10931733792, 31619704608\}
\end{displaymath}
and
\begin{displaymath}
\{-21145950, 21782754, 28598850, 56133175,
\end{displaymath}
\begin{displaymath}
386338050, 873202050, 2468426850\}
\end{displaymath}

\section{Acknowledgements}
I would like to thank the Special Collections section of the
University of Edinburgh Library for access to their copy of Gill.
I would like to thank my daughter Catriona for help with the more
difficult parts of the French texts.

\end{document}